\documentclass[a4paper,11pt,english]{smfart}

\usepackage{amssymb}
\usepackage{url}
\usepackage[T1]{fontenc}
\RequirePackage{calrsfs}
\DeclareSymbolFont{rsfscript}{OMS}{rsfs}{m}{b}
\DeclareSymbolFontAlphabet{\mathrsfs}{rsfscript}
\usepackage[utopia,expert]{mathdesign}
\usepackage{lscape}

\usepackage{enumitem}

\usepackage{todonotes}

\usepackage[leftbars]{changebar}

\usepackage{palatino}
\usepackage{rotating}
\usepackage{graphicx}

\usepackage{tikz-cd}

\usepackage{enumerate}

\usepackage{amscd}

\usepackage{color}
\definecolor{shadecolor}{gray}{0.90}

\def\proofname{D\'emonstration}

\input xypic
\xyoption{all}
\xyoption{arc}

\makeindex

\def\indexname{Index des notations}

\def\thetheo{\arabic{theo}}

\renewcommand\theequation{\thetheo}
\def\equat{\refstepcounter{theo}\begin{equation}}
\def\endequat{\end{equation}}

\renewcommand\thetable{\thetheo}

\renewcommand\thesection{\arabic{section}}
\renewcommand\thesubsection{\thesection.\Alph{subsection}}

\def\contentsname{Table des mati\`eres}
\def\listfigurename{Liste des figures}
\def\listtablename{Liste des tables}
\def\appendixname{Appendice}

\def\refname{R\'ef\'erences}

    \def\CM{{\mathbb{C}}}

    \def\FM{{\mathbb{F}}}

    \def\QM{{\mathbb{Q}}}
    
\def\SG{{\mathfrak S}}

    \def\CC{{\mathcal{C}}}

\def\Gb{{\mathbf G}}    
    \def\HC{{\mathcal{H}}}

\def\Lb{{\mathbf L}}

\def\Ob{{\mathbf O}}    
\def\Pb{{\mathbf P}}

\def\Sb{{\mathbf S}}    \def\SC{{\mathcal{S}}}
    
    \def\UC{{\mathcal{U}}}

    \def\ZC{{\mathcal{Z}}}

\def\Zrm{{\mathrm{Z}}}

\def\d{\delta}
\def\D{\Delta}
\def\e{\varepsilon}

\def\l{\lambda}
\def\L{\Lambda}

\def\O{\Omega}

\def\s{\sigma}

\def\z{\zeta}

\def\mub{{\boldsymbol{\mu}}}

\DeclareMathOperator{\Aut}{{\mathrm{Aut}}}

\def\to{\rightarrow}

\def\injto{\hookrightarrow}

\def\finl{~$\blacksquare$}

\def\lexp#1#2{\kern\scriptspace\vphantom{#2}^{#1}\kern-\scriptspace#2}
\def\le{\hspace{0.1em}\mathop{\leqslant}\nolimits\hspace{0.1em}}
\def\ge{\hspace{0.1em}\mathop{\geqslant}\nolimits\hspace{0.1em}}

\mathchardef\inferieur="321E
\mathchardef\superieur="321F

\def\eqna{\begin{eqnarray*}}
\def\endeqna{\end{eqnarray*}}

\def\itemth#1{\item[${\mathrm{(#1)}}$]}

\catcode`\@=11
\long\def\@car#1#2\@nil{#1}
\long\def\@first#1#2{#1}
\long\def\@second#1#2{#2}
\long\def\ifempty#1{\expandafter\ifx\@car#1@\@nil @\@empty
  \expandafter\@first\else\expandafter\@second\fi}
\catcode`\@=12

\theoremstyle{remark}

\theoremstyle{plain}

\def\BIL{LR}
\def\GAUCHE{L}
\def\CAR{CAR}
\def\FAM{FAM}

\def\xyinj{\ar@{^{(}->}}
\def\xysur{\ar@{->>}}

\bigskip

\makeatletter
\def\hlinewd#1{%
\noalign{\ifnum0=`}\fi\hrule \@height #1 %
\futurelet\reserved@a\@xhline}
\makeatother

\newlength\epaisLigne

\usepackage{multirow}

\makeatletter
\def\hlinewd#1{%
\noalign{\ifnum0=`}\fi\hrule \@height #1 %
\futurelet\reserved@a\@xhline}
\makeatother

\usepackage{multirow}

\usepackage{lscape}

\usepackage{pstricks}

\def\aff{{\mathrm{aff}}}
\def\sing{{\mathrm{sing}}}

\begin{document}

\title{A surface of degree 24 with 1~\!440 singularities of type $D_4$}

\author{{\sc C\'edric Bonnaf\'e}}
\address{
Institut Montpelli\'erain Alexander Grothendieck (CNRS: UMR 5149), 
Universit\'e Montpellier 2,
Case Courrier 051,
Place Eug\`ene Bataillon,
34095 MONTPELLIER Cedex,
FRANCE} 

\makeatletter
\email{cedric.bonnafe@umontpellier.fr}
%
%

\date{\today}

\thanks{The author is partly supported by the ANR 
(Project No ANR-16-CE40-0010-01 GeRepMod)}

\pagestyle{myheadings}
\markboth{\sc C. Bonnaf\'e}{Surface of degree 24 with 1~\!440 singularities of type $D_4$}

\vskip-1cm

\begin{abstract}
Using the invariant algebra of the complex reflection group denoted by $G_{32}$ 
in the Shephard--Todd classification, we construct three irreducible surfaces 
in $\Pb^3(\CM)$ with many singularities: one of them has 
degree $24$ and contains $1~\!440$ quotient singularities of type $D_4$.
\end{abstract}

\maketitle

\bigskip

Let $\mu_{D_4}(d)$ denote the maximal number of quotient singularities 
of type $D_4$ that an irreducible projective surface of degree $d$ in $\Pb^3(\CM)$ 
might have. Miyaoka~\cite{miyaoka} proved that
$$\mu_{D_4}(d) \le \frac{16}{117}d(d-1)^2.$$
For $d=8$, $16$ or $24$, this reads
$$\mu_{D_4}(8) \le 53,\qquad \mu_{D_4}(16) \le 492 \qquad\text{and}\qquad \mu_{D_4}(24) \le 1736.$$
The main results of this note are that
\equat\label{eq:8-16-24}
\mu_{D_4}(8) \ge 44,\qquad \mu_{D_4}(16) \ge 352
\qquad\text{and}\qquad \mu_{D_4}(24) \ge 1440
\endequat
and that
\equat\label{eq:8k}
\mu_{D_4}(8k) \ge 44 k^3
\endequat
for all $k \ge 1$. 
For this, let $\CM[x_1,x_2,x_3,x_4]$ be the 
polynomial ring over $\CM$ in $4$ indeterminates with its usual grading and let 
$\Pb^3(\CM) = {\mathrm{Proj}}(\CM[x_1,x_2,x_3,x_4])$ be the associated projective 
space of dimension $3$. If $f \in \CM[x_1,x_2,x_3,x_4]$ is homogeneous, 
we denote by $\ZC(f)$ the projective hypersurface it defines, and by 
$\ZC(f)_\sing$ its reduced singular locus. If $k$ is
a natural number, we denote by $f[k]$ the homogeneous polynomial 
$f(x_1^k,x_2^k,x_3^k,x_4^k)$. 

For proving the above results, we exhibit a particular homogeneous polynomial $g$ of degree $8$ such that 
the associated projective varieties $\ZC(g)$, $\ZC(g[2])$ and $\ZC(g[3])$ 
(which have respective degrees $8$, $16$ and $24$) 
have many quotient singularities of type $D_4$. The proof relies essentially 
on {\tt Magma} computations that will be detailed in the next sections: we have decided 
to write the full {\tt Magma} code in this {\tt arXiv} version, so that 
the reader can check by himself the computations, but note that this code will 
not appear in the published version of this paper.

Our polynomial $g$ is constructed using polynomial invariants of various 
finite subgroups of $\Gb\Lb_4(\CM)$. Let $W_1$ be the subgroup of $\Gb\Lb_4(\CM)$ generated by
$$
s_1=\begin{pmatrix}
0 &  1 &  0 &  0 \\
1 &  0 &  0 &  0 \\
0 &  0 &  1 &  0 \\
0 &  0 &  0 & - 1\\
\end{pmatrix},
\quad s_2 =\begin{pmatrix}
1&  0&  0&  0\\
0&  0&  1&  0\\
0&  1&  0&  0\\
0&  0&  0& - 1\\
\end{pmatrix}\quad\text{and}\quad
s_3=\begin{pmatrix}
- 1 &  0 &  0 &  0 \\
0  &  1 &  0 &  0 \\
0  &  0 &  0 &  1 \\
0  &  0 &  1 &  0 \\
\end{pmatrix}.
$$
Let $\z_3$ (resp. $\z_4$) be a primitive third (resp. fourth) root of unity. 
Let $W_2$ be the subgroup of $\Gb\Lb_4(\CM)$ generated by
$$t_1=\begin{pmatrix}
0 &  1 &  0 &  0\\
1 &  0 &  0 &  0\\
0 &  0 &  1 &  0\\
0 &  0 &  0 &  \z_4\\
\end{pmatrix},
\quad
t_2=\begin{pmatrix}
1 &  0 &  0 &  0 \\
0 &  0 &  1 &  0 \\
0 &  1 &  0 &  0 \\
0 &  0 &  0 &  \z_4 \\
\end{pmatrix}
\quad\text{and}\quad
t_3=\begin{pmatrix}
-\z_4&  0 &  0 &  0 \\
   0 &  1 &  0 &  0 \\
   0 &  0 &  0 &  1 \\
   0 &  0 &  1 &  0 \\
\end{pmatrix}.$$
Finally, let $W_3$ denote the subgroup of $\Gb\Lb_4(\CM)$ generated by 
$$u_1=\begin{pmatrix}
1 & 0 & 0    & 0 \\
0 & 1 & 0    & 0 \\
0 & 0 & \z_3 & 0 \\
0 & 0 & 0    & 1 \\
\end{pmatrix},\quad
u_2=\begin{pmatrix}
\frac{\z_3+2}{3} & \frac{\z_3-1}{3} & \frac{\z_3-1}{3} & 0 \\
\frac{\z_3-1}{3} & \frac{\z_3+2}{3} & \frac{\z_3-1}{3} & 0 \\
\frac{\z_3-1}{3} & \frac{\z_3-1}{3} & \frac{\z_3+2}{3} & 0 \\
0 & 0 & 0 & 1 \\
\end{pmatrix},$$
$$
u_3=\begin{pmatrix}
1 & 0 & 0 & 0 \\
0 & \z_3 & 0 & 0 \\
0 & 0 & 1 & 0 \\
0 & 0 & 0 & 1 \\
\end{pmatrix}
\quad\text{and}\quad
u_4=\begin{pmatrix}
\frac{\z_3+2}{3} & \frac{1-\z_3}{3} & 0 & \frac{1-\z_3}{3} \\
\frac{1-\z_3}{3} & \frac{\z_3+2}{3} & 0 & \frac{\z_3-1}{3} \\
0 & 0 & 1 & 0 \\
\frac{1-\z_3}{3} & \frac{\z_3-1}{3} & 0 & \frac{\z_3+2}{3} \\
\end{pmatrix}.$$

\medskip

\noindent{\bf Commentaries.} 
The following facts are checked using~\cite{magma}, as explained below. 
Let $\Zrm(W_i)$ denote the center of $W_i$. In all cases, 
it is isomorphic to a group of roots of unity acting by scalar 
multiplication. Then:
\begin{itemize}
\itemth{a} The group $W_1$ has order $48$ and is isomorphic to 
the non-trivial double cover $\tilde{\SG}_4$ of the symmetric group 
$\SG_4 \simeq W_1/\Zrm(W_1)$.

\itemth{b} The group $W_2$ has order $768$, 
contains a normal abelian subgroup $H$ of order $32$ 
and $W_2/H \simeq \SG_4$. The group $W_2/\Zrm(W_2)$ has order 
$192$, but is not isomorphic 
to a Coxeter group of type $D_4$.

\itemth{c} The group $W_3$ is the complex reflection group denoted 
by $G_{32}$ in the Shephard--Todd classification~\cite{sht} (it has order $155~\!920$). 
Recall that the group $W_3/\Zrm(W_3)$ is a simple group of order $25~\!920$ 
and is isomorphic to the derived subgroup of the Weyl group of type $E_6$ (i.e. 
to the derived subgroup of the special orthogonal group 
$\Sb\Ob_5(\FM_{\! 3})$). Note that we have used the form 
implemented by Michel~\cite{jean} in the {\tt Chevie} package of 
{\tt Gap3}. It contains the group $W_1$ as a subgroup, as well as 
a subgroup of diagonal matrices isomorphic to $(\mub_3)^4$, where $\mub_d$ 
is the group of $d$-th roots of unity.\finl
\end{itemize}

\bigskip

If $\l=(\l_1 \ge \l_2 \ge \l_3 \ge \l_4)$ is a partition of $8$ 
of length at most $4$, we denote by $\O_\l^-$ (resp. $\O_\l^+$) be the orbit of the monomial 
$x_1^{\l_1}x_2^{\l_2}x_3^{\l_3}x_4^{\l_4}$ under the action of $W_1$ (resp. the symemtric group 
$\SG_4$) and we set
$$m_\l^\e=\sum_{m \in \O_\l^\e} m$$
for $\e \in \{+,-\}$. 
Then $m_\l^+$ is the symmetric function traditionnally 
denoted by $m_\l$. If all the $\l_i$'s are even, then 
$m_\l^-=m_\l^+$ but note for instance that
\eqna
m_{611}^+ \neq m_{611}^- &=&  
x_1^6 x_2 x_3 + x_1^6 x_2 x_4 - x_1^6 x_3 x_4 + x_1 x_2^6 x_3 - x_1 x_2^6 x_4 + x_2^6 x_3 x_4 \\ 
&& + x_1 x_2 x_3^6 + x_1 x_3^6 x_4 - x_2 x_3^6 x_4 - x_1 x_2 x_4^6 - x_1 x_3 x_4^6 - x_2 x_3 x_4^6. 
\endeqna
Now, let
\eqna
g\!\!\!\!&=&\!\!\!\!m_8^- - 6 m_{62}^- - 60 m_{611}^- + 2~\!240 m_{521}^- -14 m_{44}^- 
+ 10~\!180 m_{431}^- + 40~\!412 m_{422}^- \\ && - 23~\!440 m_{4211}^- 
+ 111~\!980 m_{332}^- + 154~\!704 m_{2222}^-.
\endeqna
By construction, $m_\l^-$ is invariant under the action of $W_1$ and so 
$g$ is invariant under the action of $W_1 \simeq \tilde{\SG}_4$. One can check with {\tt Magma} 
the following facts:

\bigskip

\noindent{\bf Proposition 1.} {\it If $1 \le k \le 3$, then the polynomial $g[k]$ 
is invariant under the action of $W_k$.}

\bigskip

\noindent{\bf Theorem 2.} {\it The homogeneous polynomial $g$ satisfies the following statements:
\begin{itemize}
\itemth{a} $\ZC(g)$ is an irreducible surface of degree $8$ in $\Pb^3(\CM)$ 
with exactly $44$ singular points which are all quotient singularities of 
type $D_4$.

\smallskip

\itemth{b} If $k \ge 1$, then $\ZC(g[k])$ is an irreducible surface 
of degree $8k$, whose singular locus has dimension $0$ and contains at least $44k^3$ quotient 
singularities of type $D_4$.

\smallskip

\itemth{c} $\ZC(g[2])$ is an irreducible surface of degree $16$ with exactly $472$ singular points: 
$24$ quotient singularities of type $A_1$, $96$ quotient singularities of type $A_2$ and 
$352$ quotient singularities of type $D_4$.

\smallskip

\itemth{d} $\ZC(g[3])$ is an irreducible surface of degree $24$ in $\Pb^3(\CM)$ 
with exactly $1~\!440$ singular points which are all quotient singularities of type $D_4$. 
The automorphism group of $\ZC(g[3])$ contains at least $25~\!920$ elements 
and acts transitively on the $1~\!440$ singular points.
\end{itemize}}

\bigskip

\noindent{\bf Remark 1.} 
It turns out that we did not find the polynomial $g$ directly: 
we found first $g[3]$ by looking at invariants of degree $24$ of $W_3 \simeq G_{32}$, 
following ideas of Barth~\cite{barth} and Sarti~\cite{sarti 0},~\cite{sarti 1},~\cite{sarti} 
(for constructing the {\it Barth sextic} with $65$ nodes and, for instance, 
the {\it Sarti dodecic} with $600$ nodes), 
who used invariants of Coxeter groups of type $H_3$ and $H_4$. 
See also~\cite{bonnafe curves} for details about the method used for finding $g$.\finl

\bigskip

\noindent{\bf Remark 2.} Note that $g$ has coefficients in $\QM$ but the singular points 
of $\ZC(g)$, $\ZC(g[2])$ and $\ZC(g[3])$ have coordinates 
in various field extensions of $\QM$, and most of the singular points 
are not real (at least in this model).\finl

\bigskip

So let us start by defining the polynomial $g$ and the three groups 
$W_i$ in {\sc Magma} and checking the facts stated in Proposition~1 
and Commentaries. These data, together with the definition of the fields 
$K$ and $L$ as well as a function {\tt ProjectiveOrbit} 
used for computing orbits of various points under the action ot the groups $W_i$, 
are all contained in a file {\tt g32-article.m}, whose content is given 
in the Appendix. Note that we will first work with the projective space over 
$\QM$, and the polynomial $g$ will be defined over $\QM$.

\medskip

\begin{quotation}
\small
\begin{verbatim}
> load 'g32-article.m';
> 
> Order(W1);
48
> Order(Centre(W1));
2
> bool:=IsIsomorphic(W1/Centre(W1),SymmetricGroup(4));
> bool;
true
> Order(DerivedSubgroup(W1));
24
> gK:=ChangeRing(g,K);
> [gK^i eq gK : i in Generators(W1)];
[ true, true, true ]
> 
> Order(W2);
768
> Order(Centre(W2));
4
> g2L:=ChangeRing(g2,L);
> [g2L^i eq g2L : i in Generators(W2)];
[ true, true, true ]
>  
> Order(W3);
155520
> Order(W3/Centre(W3));
25920
> IsSimple(W3/Centre(W3));
true
> g3K:=ChangeRing(g3,K);
> [g3K^i eq g3K : i in Generators(W3)];
[ true, true, true, true ]
\end{verbatim}
\end{quotation}

\medskip

We now turn to the study of the singularities 
of the varieties $\ZC(g[i])$ for $i \in \{1,2,3\}$. Note the 
following fact, that will be used further:

\bigskip

\noindent{\bf Lemma 3.} {\it If $1 \le i < j \le 4$, then the closed 
subscheme of $\Pb^3(\CM)$ defined by the ideal 
$\langle g,\frac{\partial g}{\partial x_i},\frac{\partial g}{\partial x_j} \rangle$ 
has dimension $0$.}

\bigskip

This is checked thanks to the following code:

\medskip

\begin{quotation}
\small
\begin{verbatim}
> dg:=[Derivative(g,i) : i in [1..4]];
> pairs:=[[1,2],[1,3],[1,4],[2,3],[2,4],[3,4]];
> time [Dimension(Scheme(P3,[g,dg[i[1]],dg[i[2]]])) : 
>   i in pairs];
[ 0, 0, 0, 0, 0, 0 ]
Time: 20.510  // about 20 seconds
\end{verbatim}
\end{quotation}

\bigskip

\section{Degree ${\boldsymbol{8}}$}\label{sec:8}

\medskip

The {\tt Magma} computations leading to the 
proof of the statement~(a) of Theorem~2 are detailed in this section. 
Along these computations, the following facts are obtained 
(here, $\UC$ denotes the open subset of $\Pb^3(\CM)$ defined by $x_1x_2x_3x_4 \neq 0$):

\bigskip

\noindent{\bf Proposition 4.} {\it We have:
\begin{itemize}
\itemth{a} $\dim \ZC(g)_\sing = 0$, so $\ZC(g)$ is irreducible. 

\itemth{b} $\ZC(g)_\sing$ is contained in $\UC$.

\itemth{c} The group $W_1$ has $3$ orbits in $\ZC(g)_\sing$, of respective length $8$, $12$ and $24$. 
\end{itemize}}

\bigskip

We first check that $\ZC(g)_\sing$ has dimension $0$ and is contained 
in $\UC$. 

\medskip

\begin{quotation}
\small
\begin{verbatim}
> Zg:=Surface(P3,g);
> Zgsing:=SingularSubscheme(Zg);
> Dimension(Zgsing);
0
> H:=Scheme(P3,x1*x2*x3*x4);
> Dimension(Intersection(Zgsing,H));
-1
\end{verbatim}
\end{quotation}

\medskip

\noindent In particular, all the singular points 
are contained (for instance) in the affine chart $\ZC(g)^\aff$ defined by ``$x_4 \neq 0$''. 
We will make all the remaining computations in this affine chart (and extend the scalars to the field $K$):

\medskip

\begin{quotation}
\small
\begin{verbatim}
> Zgaff:=AffinePatch(Zg,1);
> ZgaffK:=ChangeRing(Zgaff,K);
> ZgaffKsing:=SingularSubscheme(ZgaffK);
> irr1:=IrreducibleComponents(ZgaffKsing);
> irr1:=[ReducedSubscheme(i) : i in irr1];
> Set([Degree(i) : i in irr1]);
{ 1 }
> ZgK:=ChangeRing(Zg,K);
> sings1:=[[Coordinates(i) : i in RationalPoints(j)] : 
>   j in irr1];
> sings1:=&cat sings1;
> sings1:=[ZgK ! (i cat [1]) : i in sings1];
> # sings1;
44
\end{verbatim}
\end{quotation}

\medskip

\noindent The result of the command {\tt Set([Degree(i) : i in irrg])} 
shows that all singular points have coordinate 
in $K$, and the last command shows that the number of singular points in $\ZC(g)$ is equal to 
$44$. 
We then determine the $W_1$-orbits in $\ZC(g)_\sing$ and check 
that they are all quotient singularities of type $D_4$ by 
picking up one point in each orbit.

\medskip

\begin{quotation}
\small
\begin{verbatim}
> orbits:=[];
> test:=sings1;
> while (# test) gt 0 do
while>   orb:=ProjectiveOrbit(W1,test[1]);
while>   orb:=[ZgK ! Coordinates(i) : i in orb];
while>   orbits:=orbits cat [orb];
while>   test:=[i : i in test | (i in orb) eq false];
while> end while;
> [# i : i in orbits]; 
[ 24, 8, 12 ]
> for i in orbits do
for>   print IsSimpleSurfaceSingularity(i[1]);
for> end for;
true D 4
true D 4
true D 4
\end{verbatim}
\end{quotation}

\medskip

\begin{center}
\begin{figure}
\begin{center}
\includegraphics[scale=0.2]{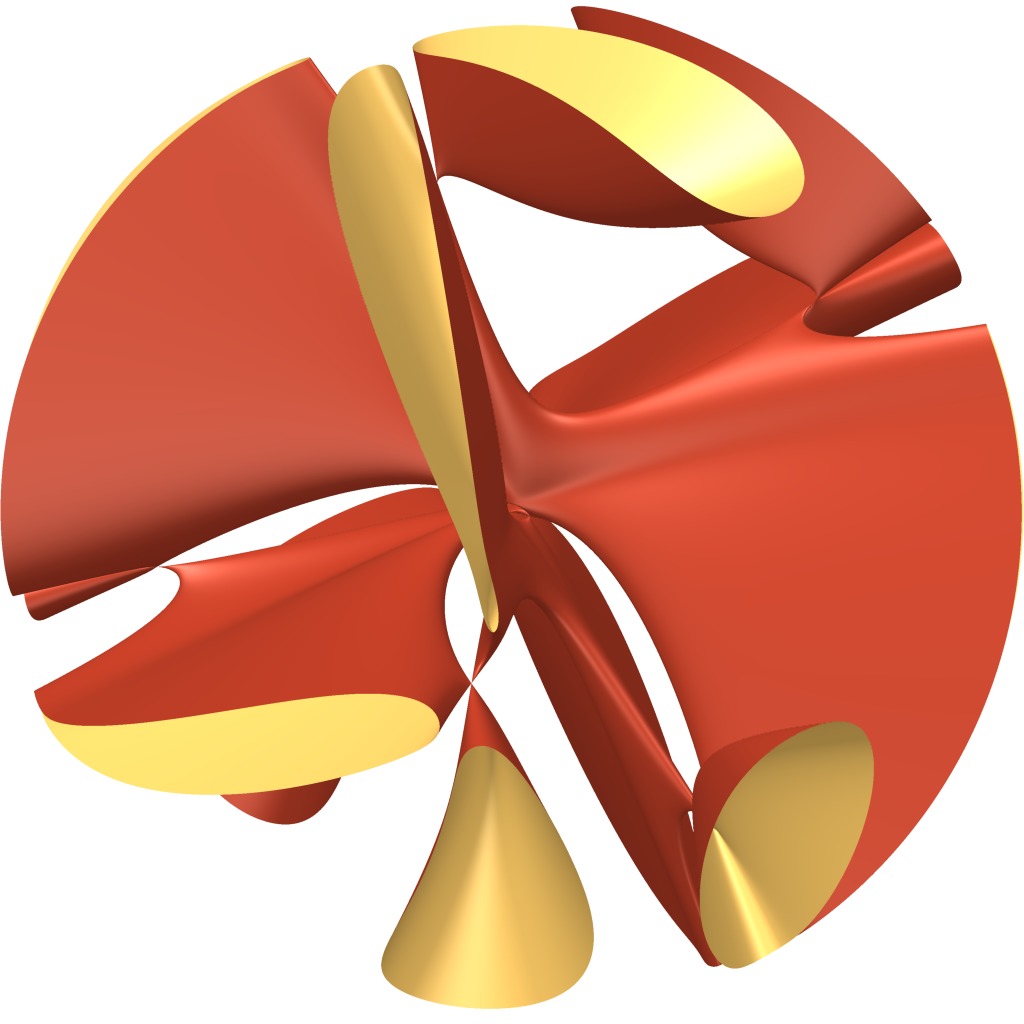}
\includegraphics[scale=0.2]{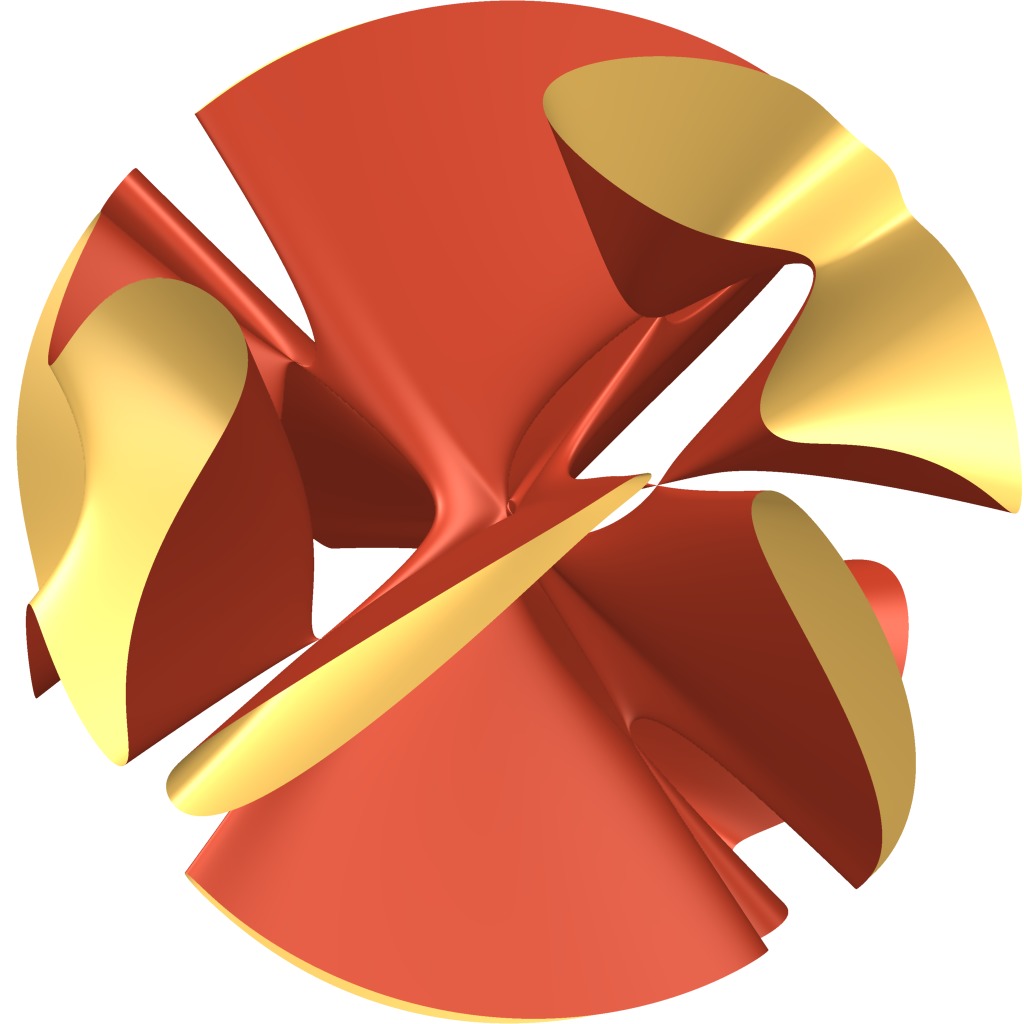}
\caption{Part of the real locus of $\ZC(g)$}\label{fig:8}
\end{center}
\end{figure}
\end{center}

Note that the points in the $W_1$-orbit of cardinality $8$ 
are the only real singular points of $\ZC(g)$. 
Figure~\ref{fig:8} shows part of the real locus of $\ZC(g)$. 

\bigskip

\section{Degree ${\boldsymbol{8k}}$}\label{sec:8k}

\medskip

Let $\UC$ denote the open subset of $\Pb^3(\CM)$ defined by $x_1x_2x_3x_4 \neq 0$ 
and let $\s_{\! k} : \Pb^3(\CM) \to \Pb^3(\CM)$, $[x_1;x_2;x_3;x_4] \mapsto [x_1^k;x_2^k;x_3^k;x_4^k]$. 
The restriction of $\s_{\! k}$ to a morphism $\UC \to \UC$ is an \'etale Galois covering, 
with group $(\mub_k)^4/\D\mub_k$ (here, $\D : \mub_k \injto (\mub_k)^4$ 
is the diagonal embedding). We have $\ZC(g[k])=\s_{\! k}^{-1} (\ZC(g))$. 

Let us first prove that $\ZC(g[k])$ is irreducible. We may assume that 
$k \ge 2$, as the result has been proved for $k=1$ in the previous section. 
Recall that
$$\frac{\partial g[k]}{\partial x_i}=kx_i^{k-1}(\frac{\partial g}{\partial x_i} \circ \s_{\! k}),$$
so the singular locus of $\ZC(g[k])$ is contained in 
$$\{p_1,p_2,p_3,p_4\} \cup \Bigl(\bigcup_{i \neq j} \s_{\! k}^{-1} (\ZC_{i,j})\Bigr),$$
where $p_i=[\d_{i1};\d_{i2};\d_{i3};\d_{i4}]$ (and $\d_{ij}$ is the Kronecker symbol) and 
$\ZC_{i,j}$ is the subscheme of $\Pb^3(\CM)$ defined by the ideal 
$\langle g,\frac{\partial g}{\partial x_i},\frac{\partial g}{\partial x_j} \rangle$ 
(and which has dimension $0$ by Lemma~3). 
Since $\s_{\! k}$ is finite, this implies that $\ZC(g[k])_\sing$ has dimension $0$, 
so $\ZC(g[k])$ is irreducible.

Now, $\s_{\! k} : \UC \to \UC$ is \'etale and the singular locus of $\ZC(g)$ 
is contained in $\UC$ (see Proposition~4(b)). Therefore, the $44$ singularities of $\ZC(g)$ lift to  
$44 k^3$ singularities in $\ZC(g[k]) \cap \UC$ of the same type, i.e. quotient 
singularities of type $D_4$. This proves the statement~(b) of Theorem~2.

Note that, for $k=2$, $3$ and $4$ (and maybe for bigger $k$) we will prove 
in the next sections that $\ZC(g[k])$ contains singular points outside 
of $\UC$. 

\bigskip

\section{Degree ${\boldsymbol{16}}$}\label{sec:16}

\medskip

Using the morphism $\s_{\! 2}$ defined in the previous section, we get that $\ZC(g[2]) \cap \UC$ 
has exactly $352$ singular points, which are all quotient singularities 
of type $D_4$. We now need to determine the singularities which are not contained 
in $\UC$. So let $\HC$ be the complement of $\UC$ in $\Pb^3(\CM)$. 

\medskip

\begin{quotation}
\small
\begin{verbatim}
> Zg2:=Surface(P3,g2); 
> Zg2sing:=SingularSubscheme(Zg2);
> H:=Scheme(P3,x1*x2*x3*x4);
> Zg2singH:=Intersection(Zg2sing,H);
> time Zg2singH:=ReducedSubscheme(Zg2singH);
Time: 11.300
> irr2H:=IrreducibleComponents(Zg2singH);
> # irr2H;
18
> &+ [Degree(i) : i in irr2H];
120
\end{verbatim}
\end{quotation}

\medskip

\noindent The last command shows that $\ZC(g)_\sing \cap \HC$ contains $120$ points. 
We now check that all the singular points contained in $\HC$ have coordinates in $L$:

\medskip

\begin{quotation}
\small
\begin{verbatim}
> Zg2L:=ChangeRing(Zg2,L);
> sings2H:=[[Coordinates(i) : 
>   i in RationalPoints(ChangeRing(j,L))] : 
>   j in irr2H];
> sings2H:=&cat sings2H;
> sings2H:=[Zg2L ! i : i in sings2H];
> # sings2H;
120
\end{verbatim}
\end{quotation}

\medskip

\noindent We now determine the $W_2$-orbits in $\ZC(g)_\sing \cap \HC$: 
there is one $W_2$-orbit of cardinality $24$ (and we check that 
its elements are quotient singularities of type $A_1$) and one
of cardinality $96$. 

\medskip

\begin{quotation}
\small
\begin{verbatim}
> orbits2:=[];
> test:=sings2H;
> W2L:=ChangeRing(W2,L);
> while (# test) gt 0 do
while>   orb:=ProjectiveOrbit(W2L,test[1]);
while>   orb:=[Zg2L ! Coordinates(i) : i in orb];
while>   orbits2:=orbits2 cat [orb];
while>   test:=[i : i in test | (i in orb) eq false];
while> end while;
> [# i : i in orbits2];
[ 24, 96 ]
> IsSimpleSurfaceSingularity(orbits2[1][1]);
true A 1
\end{verbatim}
\end{quotation}

\medskip

\noindent We now study the singularity of $\ZC(g[2])$ at the points of the orbit 
of cardinality $96$. It turns out that that the command {\tt IsSimpleSurfaceSingularity} 
takes too much time to get a conclusion, so we will investigate properties 
of the equation of $\ZC(g[2])$ in a neighborhood of the first point $p$ 
(in {\tt Magma} list {\tt orbits[2]}). We work in the affine chart $x_3 \neq 0$ 
(where $p=[\xi_1;\xi_2;\xi_3;\xi_4]$ lives), and we denote by $(x,y,z)$ 
the coordinates of the affine chart $x_3 \neq 0$ equal to 
$(x_1/x_3+\xi_1/\xi_3,x_2/x_3+\xi_2/\xi_3,x_4/x_3+\xi_4/\xi_3)$ and we set
$$f(x,y,z)=g[2](x,y,1,z).$$
If $j \ge 0$, we denote by $f_j$ the homogeneous component of $f$ of degree $j$. 
As $p \in \ZC(g[2])_\sing$, we have $f_0=f_1=0$.

\medskip

\begin{quotation}
\small
\begin{verbatim}
> p:=orbits2[2][1];
> A3L<x,y,z>:=AffineSpace(L,3);
> cop:=Coordinates(p);
> f:=Evaluate(g2L,[x+cop[1],y+cop[2],1,z+cop[4]]);
> cof:=Coefficients(f);
> mof:=Monomials(f);
> l:=# mof;
> f2:=&+ [cof[i]*mof[i] : i in [1..l] | 
>   Degree(mof[i]) eq 2];
> Factorization(f2);
[
    <y^2 + 1/4550725*(-298386*xi^7 + 4375808*xi^6
    - 795576*xi^5 + 3978000*xi^4 - 99422*xi^3 
    - 1989000*xi^2 + 696114*xi - 3679654)*alpha*y*z 
    + 1/182029*(201619*xi^6 - 403238*xi^2 - 472435)*z^2, 1>
]
\end{verbatim}
\end{quotation}

\medskip

\noindent The last command 
shows that there exists a linear change of the coordinates $(x,y,z) \mapsto (x,Y,Z)$ 
such that $f_2$ might be transformed into $Y^2+Z^2$. By standard arguments, this proves that $p$ 
is a quotient singularity of type $A_k$, for some $k \ge 2$, 
which can be obtained as the Milnor number of $f$: 
however, {\tt Magma} cannot compute this Milnor number 
in a reasonable amount of time and we need to copy 
the polynomial $f$ in the software {\tt Singular}~\cite{singular} 
to compute this Milnor number (!): we obtain $2$. So $p$ is a quotient singularity 
of type $A_2$. This completes the proof of statement~(c) of Theorem~2.

\begin{center}
\begin{figure}
\begin{center}
\includegraphics[scale=0.2]{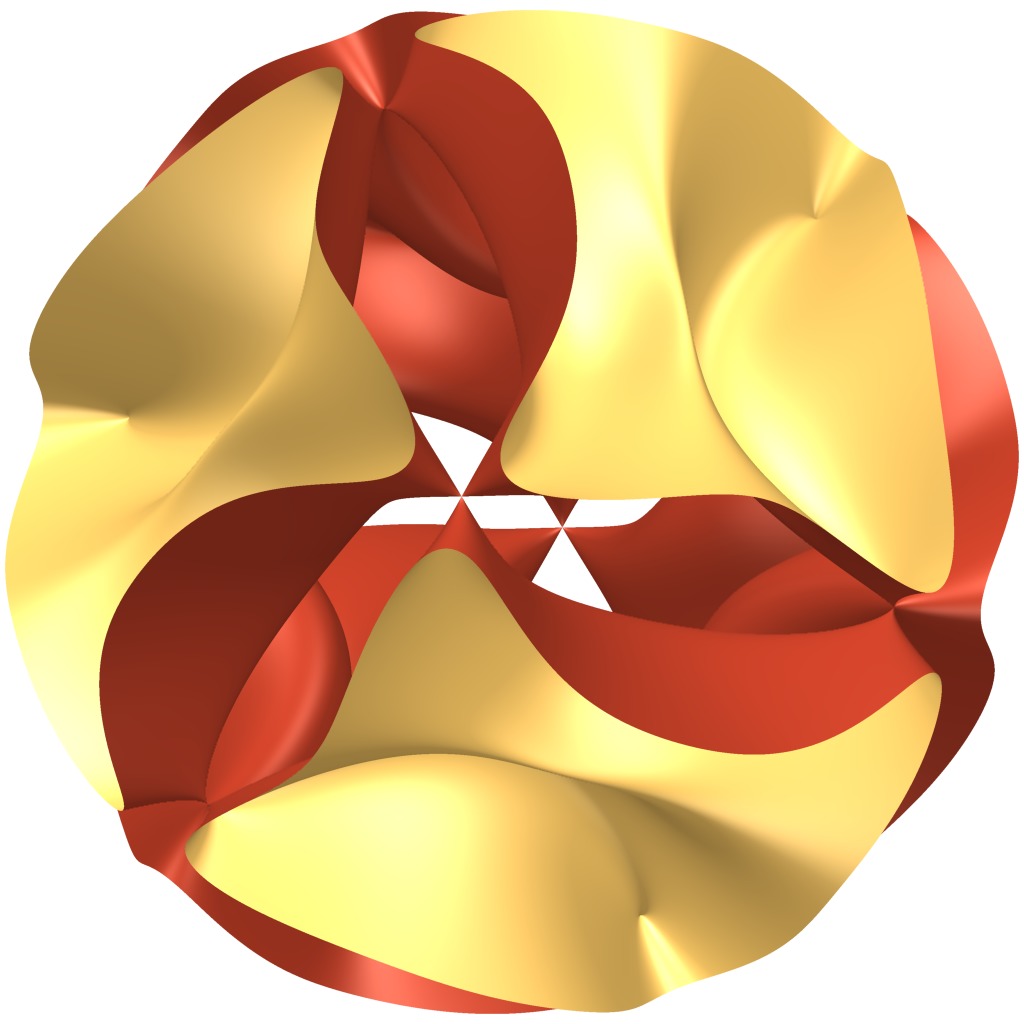}
\includegraphics[scale=0.2]{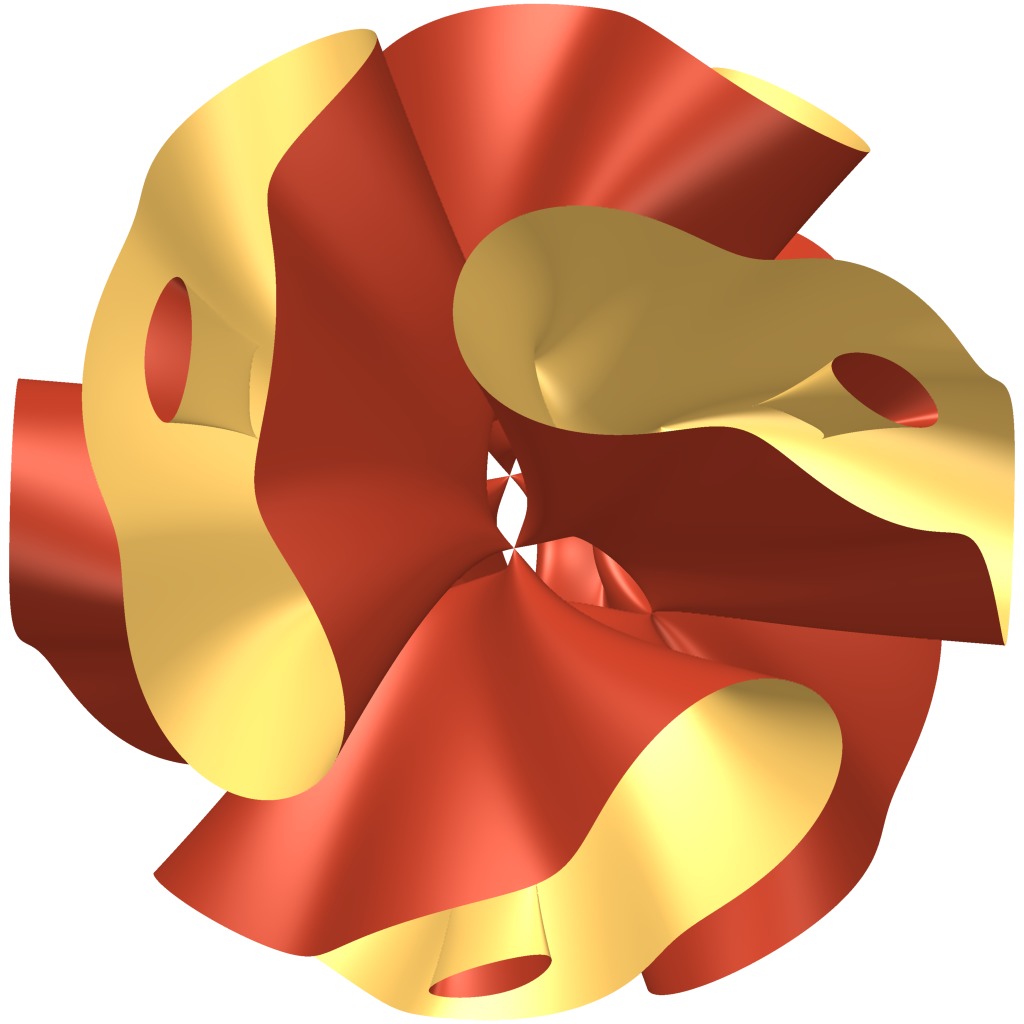}
\caption{Part of the real locus of $\ZC(g[2])$}\label{fig:16}
\end{center}
\end{figure}
\end{center}

Figure~\ref{fig:16} shows part of the real locus of $\ZC(g)$.

\section{Degree ${\boldsymbol{24}}$}\label{sec:24}

\medskip

Using the morphism $\s_{\! 3}$ defined in Section~\ref{sec:8k}, we get that $\ZC(g[3]) \cap \UC$ 
has exactly $44 \times 3^3=1~\!188$ singular points, which are all quotient singularities 
of type $D_4$. Let us 
compute $\ZC(g[3])_\sing \cap \HC$:

\medskip

\begin{quotation}
\small
\begin{verbatim}
> Zg3:=Surface(P3,g3); 
> Zg3sing:=SingularSubscheme(Zg3);
> Zg3singH:=Intersection(Zg3sing,H);
> time Zg3singH:=ReducedSubscheme(Zg3singH);
Time: 19.320
> time irr3H:=IrreducibleComponents(Zg3singH);
Time: 18.170
> # irr3H;
72
> Set([Degree(i) : i in irr3H]);
{ 2, 4 }
> &+ [Degree(i) : i in irr3H];
252
\end{verbatim}
\end{quotation}

\medskip

\noindent The last command shows that $\ZC(g[3])_\sing \cap \HC$ contains 
$252$ points. We now show that they are all defined over $K$:

\medskip

\begin{quotation}
\small
\begin{verbatim}
> Zg3K:=ChangeRing(Zg3,K);
> sings3H:=[[Coordinates(i) : 
>   i in RationalPoints(ChangeRing(j,K))] : 
>   j in irr3H];
> sings3H:=&cat sings3H;
> sings3H:=[Zg3K ! i : i in sings3H];
> # sings3H;
252
\end{verbatim}
\end{quotation}

\medskip

\noindent So $\ZC(g[3])_\sing$ contains $1440$ points, and we now check that 
$W_3$ acts transitively on them:

\medskip

\begin{quotation}
\small
\begin{verbatim}
> p:=sings3H[1];
> time # ProjectiveOrbit(W3,p);
1440
Time: 29.850
\end{verbatim}
\end{quotation}

\medskip

\begin{figure}
\includegraphics[scale=0.2]{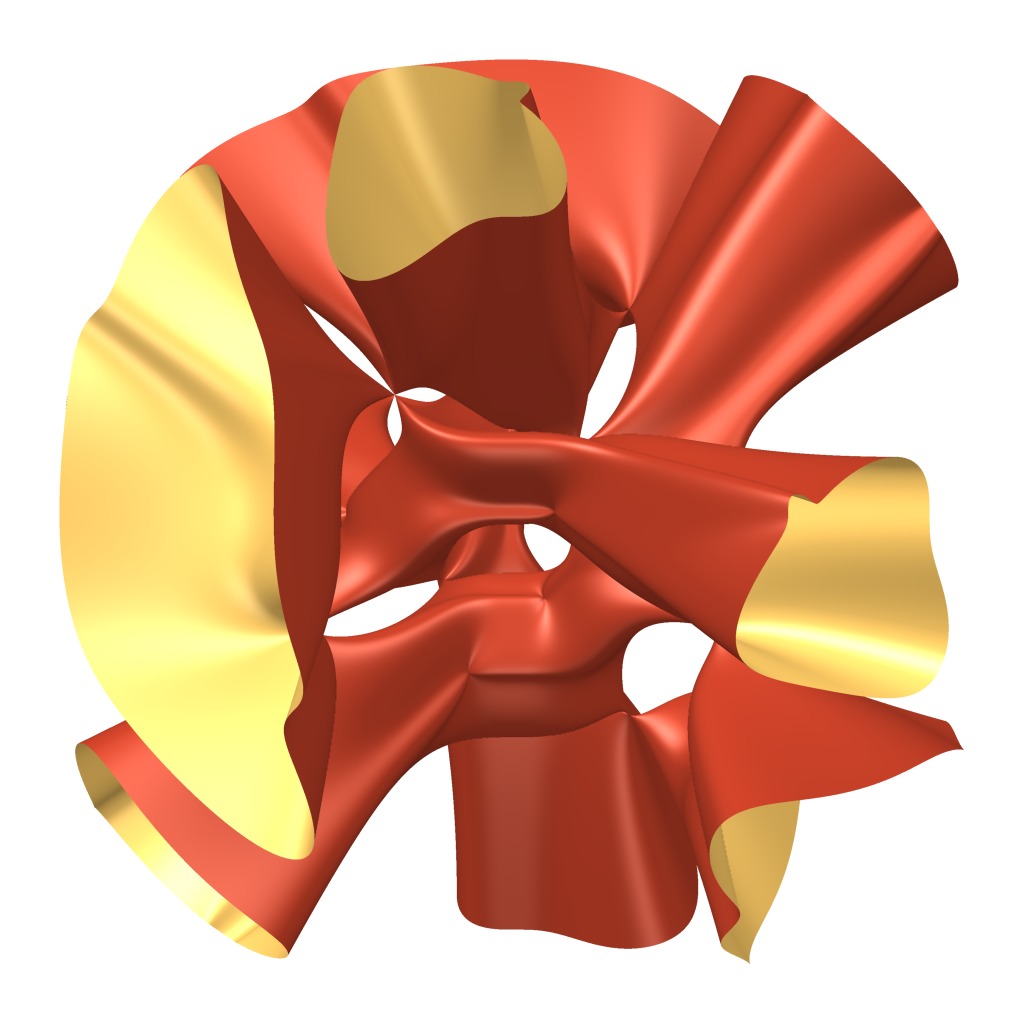}
\includegraphics[scale=0.2]{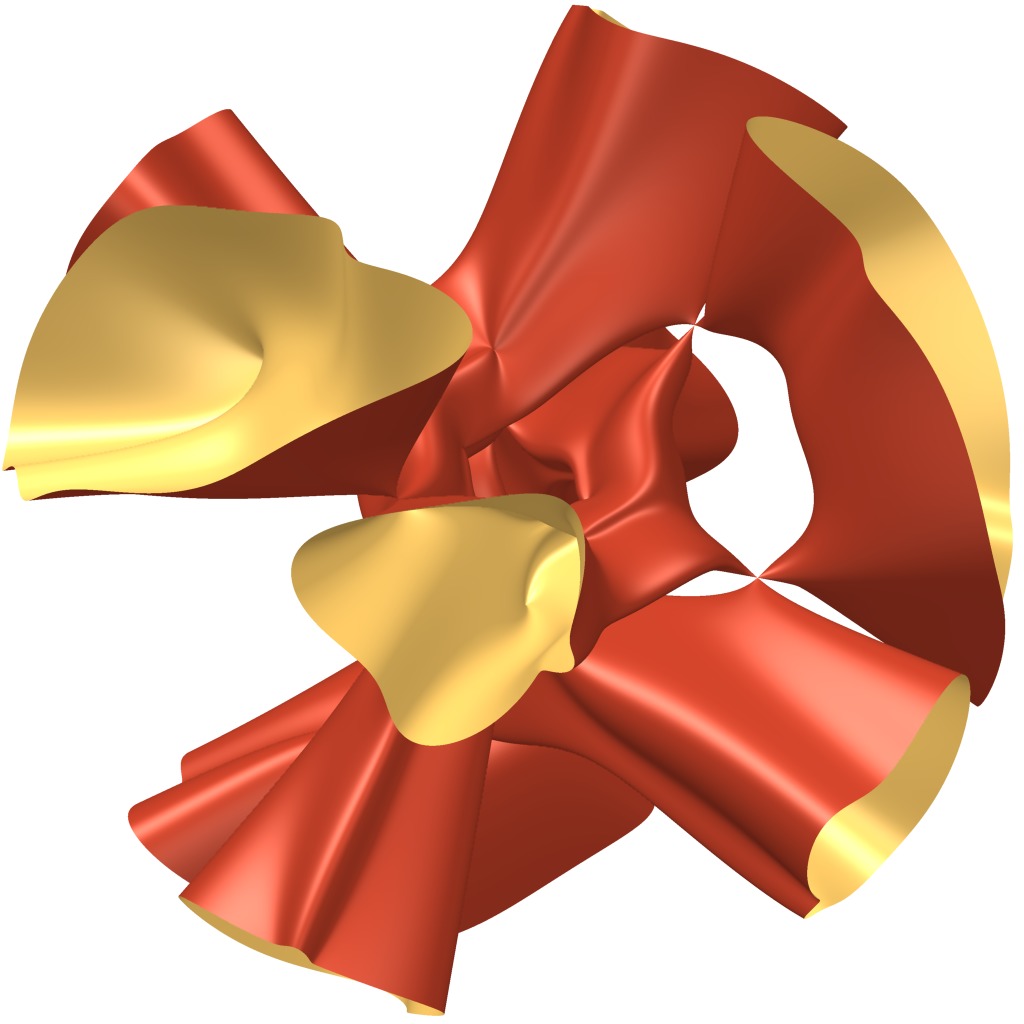}
\caption{Part of the real locus of $\ZC(g[3])$}\label{fig:24}
\end{figure}

\noindent The proof 
of statement~(d) of Theorem~2 is complete: Figure~\ref{fig:24} gives partial views 
of its real locus.

\bigskip

\section{Complements}

\medskip

\noindent{\bf Remark 3.} 
From Section~\ref{sec:8k}, we deduce that $\ZC(g[4])_\sing$ has $2~\!816$ quotient 
singularities of type $D_4$ lying in the open subset $\UC$ and it can be checked that 
it has $480$ other singular points not lying in $\UC$, for which we did not determine
the type.\finl

\bigskip

\noindent{\bf Remark 4.} 
After investigations in the invariant rings of several irreducible primitive complex reflection groups 
(there are $34$ such groups, denoted by $G_i$ with $4 \le i \le 37$ in Shephard--Todd classification~\cite{sht}), 
we have also been able to construct curves with many singularities. 
For example:
\begin{itemize}
\item[$\bullet$] Using the reflection group $W=G_{24}$, we have obtained 
a cuspidal curve $\CC_{14}$ of degree $14$ in $\Pb^2(\CM)$ 
with exactly $42$ cusps (all lying in a single $W$-orbit). Note that 
$W/\Zrm(W) \subset \Aut(\CC_{14})$ has order $168$ and is isomorphic to $\Gb\Lb_3(\FM_{\! 2})$.

\item[$\bullet$] Using the reflection group $W=G_{26}$, we have obtained 
a curve $\CC_{18}$ of degree $18$ in $\Pb^2(\CM)$ 
with $72$ cusps and $12$ nodes (these are the two $W$-orbits of singular points). 
Note that $W/\Zrm(W) \subset \Aut(\CC_{18})$ has order $216$.
\end{itemize}
Also, other singular surfaces have been obtained. 
For example:
\begin{itemize}
\item[$\bullet$] Using the reflection group $W=G_{29}$ 
(note that $W/\Zrm(W)$ has order $1~\!920$), we have obtained:
\begin{itemize}
\item[-] a surface of degree $8$ in $\Pb^3(\CM)$ with $160$ nodes, 
all belonging to the same $W$-orbit. Recall that the {\it Endra\ss~octic}~\cite{endrass}  
has degree $8$ and $168$ nodes (and its automorphism group 
has order $16$) and the {\it Sarti octic}~\cite{sarti 0} has $144$ nodes.

\item[-] a surface of degree $8$ in $\Pb^3(\CM)$ with $20$ singular points of multiplicity 
$3$ and Milnor number $11$, all belonging to the same $W$-orbit. 
\end{itemize}

\item[$\bullet$] Using the reflection group $W=G_{31}$, 
we have obtained a surface $\SC_{20}$ of degree $20$ in $\Pb^3(\CM)$ with $1~\!920$ nodes, all lying in 
the same $W$-orbit. Note that $W/\Zrm(W) \subset \Aut(\SC_{20})$ 
has order $11~\!520$) and recall that the {\it Chmutov surface}~\cite{chmutov} of 
degree $20$ has $2~\!926$ nodes.
\end{itemize}
Details will appear in a forthcoming paper~\cite{bonnafe curves}.\finl

\bigskip

\noindent{\bf Acknowledgements.} This paper is based 
upon work supported by the National Science Foundation under 
Grant No. DMS-1440140 while the author was in residence at the Mathematical Sciences Research 
Institute in Berkeley, California, during the Spring 2018 semester. 
The hidden computations which led to the discovery 
of the polynomial $g$ were done using the High Performance Computing facilities 
of the MSRI.

I wish to thank warmly Alessandra Sarti, Oliver Labs and Duco van Straten 
for useful comments and references and Gunter Malle for a careful reading 
of a first version of this note. 
Figures were realized using the software {\sc SURFER}~\cite{surfer}.

\bigskip

\section*{Appendix}

\medskip

This appendix gives a copy of the file {\tt g32-article.m} loaded at the beginning 
of the computations. It contains the data of the polynomial $g$, the fields 
$K$ and $L$, the three groups 
$W_i$ and the function {\tt ProjectiveOrbit} which is used throughout the 
computations (it is certainly not the most efficient code, but it is 
sufficient for our purpose). Note that 
$W_1$ and $W_3$ are defined over the field $K$, 
while $W_2$ is defined over the field $L$.

\medskip

\begin{quotation}
\small
\begin{verbatim}
Q:=RationalField();
P3<x1,x2,x3,x4>:=ProjectiveSpace(Q,3);
K<zeta>:=CyclotomicField(12);
zeta3:=zeta^4;
P3K:=ProjectiveSpace(K,3);
K24<xi>:=CyclotomicField(24);
zeta4:=xi^6;
POL<T>:=PolynomialRing(K24);
L<alpha>:=NumberField(T^2-(18*xi^6 + 14*xi^5 
  + 48*xi^4 + 2*xi^3 - 36*xi^2 - 14*xi - 24));
P3L:=ProjectiveSpace(L,3);
\end{verbatim}
\end{quotation}

\newpage

\begin{quotation}
\small
\begin{verbatim}
g:=x1^8 + x2^8 + x3^8 + x4^8 
- 6*(x1^6*x2^2 + x1^6*x3^2 + x1^6*x4^2 + x1^2*x2^6 
   + x2^6*x3^2 + x2^6*x4^2 + x1^2*x3^6 + x2^2*x3^6 
   + x3^6*x4^2 + x1^2*x4^6 + x2^2*x4^6 + x3^2*x4^6)
- 60*(x1^6*x2*x3 + x1^6*x2*x4 - x1^6*x3*x4 
    + x1*x2^6*x3 - x1*x2^6*x4 + x2^6*x3*x4 
    + x1*x2*x3^6 + x1*x3^6*x4 - x2*x3^6*x4
    - x1*x2*x4^6 - x1*x3*x4^6 - x2*x3*x4^6)
+ 2240*(x1^5*x2^2*x3 - x1^5*x2^2*x4 + x1^5*x2*x3^2 
      - x1^5*x2*x4^2 + x1^5*x3^2*x4 - x1^5*x3*x4^2 
      + x1^2*x2^5*x3 + x1^2*x2^5*x4 + x1*x2^5*x3^2 
      + x1^2*x2*x3^5 - x1*x2^5*x4^2 + x2^2*x3*x4^5 
      - x2*x3^5*x4^2 + x2^2*x3^5*x4 + x1^2*x2*x4^5 
      - x2*x3^2*x4^5 - x1*x2^2*x4^5 - x1^2*x3*x4^5 
      - x1^2*x3^5*x4 + x1*x2^2*x3^5 - x1*x3^5*x4^2
      + x1*x3^2*x4^5 - x2^5*x3^2*x4 - x2^5*x3*x4^2)
- 14*(x1^4*x2^4 + x1^4*x3^4 + x1^4*x4^4 
    + x2^4*x3^4 + x2^4*x4^4 + x3^4*x4^4)
+ 10180*(x1^4*x2^3*x3 + x1^4*x2^3*x4 + x1^4*x2*x3^3 
       + x1^4*x2*x4^3 - x1^4*x3^3*x4 - x1^4*x3*x4^3 
       + x1^3*x2^4*x3 - x1^3*x2^4*x4 + x1^3*x2*x3^4 
       - x1^3*x2*x4^4 + x1^3*x3^4*x4 + x1*x2^4*x3^3
       + x1*x3^4*x4^3 - x1*x3^3*x4^4 + x2^4*x3*x4^3 
       - x2^3*x3^4*x4 - x1*x2^4*x4^3 + x1*x2^3*x3^4 
       - x1*x2^3*x4^4 - x2^3*x3*x4^4 - x2*x3^4*x4^3 
       - x2*x3^3*x4^4 + x2^4*x3^3*x4 - x1^3*x3*x4^4 )
+ 40412*(x1^4*x2^2*x3^2 + x1^4*x2^2*x4^2 + x1^4*x3^2*x4^2 
       + x1^2*x2^4*x3^2 + x1^2*x2^4*x4^2 + x1^2*x3^4*x4^2 
       + x1^2*x3^2*x4^4 + x1^2*x2^2*x3^4 + x1^2*x2^2*x4^4 
       + x2^4*x3^2*x4^2 + x2^2*x3^4*x4^2 + x2^2*x3^2*x4^4)
- 23440*(x1^4*x2^2*x3*x4 - x1^4*x2*x3^2*x4 
       - x1^4*x2*x3*x4^2 + x1^2*x2*x3^4*x4 
       + x1^2*x2*x3*x4^4 - x1^2*x2^4*x3*x4 
       - x1*x2^4*x3*x4^2 - x1*x2^2*x3^4*x4 
       + x1*x2^2*x3*x4^4 - x1*x2*x3^4*x4^2 
       + x1*x2*x3^2*x4^4 + x1*x2^4*x3^2*x4 )
+ 111980*(x1^3*x2^3*x3^2 - x1^3*x2^3*x4^2 + x1^3*x2^2*x3^3 
        - x1^3*x2^2*x4^3 - x1^3*x3^3*x4^2 + x1^3*x3^2*x4^3 
        + x1^2*x2^3*x3^3 + x1^2*x2^3*x4^3 - x1^2*x3^3*x4^3 
        - x2^3*x3^3*x4^2 - x2^3*x3^2*x4^3 + x2^2*x3^3*x4^3) 
+ 154704*x1^2*x2^2*x3^2*x4^2;
g2:=Evaluate(g,[x1^2,x2^2,x3^2,x4^2]);
g3:=Evaluate(g,[x1^3,x2^3,x3^3,x4^3]);
\end{verbatim}
\end{quotation}

\newpage

\begin{quotation}
\small
\begin{verbatim}
s1:=Matrix(K,4,4,
[[ 0,  1,  0,  0],
 [ 1,  0,  0,  0],
 [ 0,  0,  1,  0],
 [ 0,  0,  0, -1]]);
s2:=Matrix(K,4,4,
[[ 1,  0,  0,  0],
 [ 0,  0,  1,  0],
 [ 0,  1,  0,  0],
 [ 0,  0,  0, -1]]);
s3:=Matrix(K,4,4,
[[-1,  0,  0,  0],
 [ 0,  1,  0,  0],
 [ 0,  0,  0,  1],
 [ 0,  0,  1,  0]]);
W1:=MatrixGroup<4,K | [s1,s2,s3]>;

t1:=Matrix(L,4,4,
[[ 0,  1,  0,     0],
 [ 1,  0,  0,     0],
 [ 0,  0,  1,     0],
 [ 0,  0,  0, zeta4]]);
t2:=Matrix(L,4,4,
[[ 1,  0,  0,     0],
 [ 0,  0,  1,     0],
 [ 0,  1,  0,     0],
 [ 0,  0,  0, zeta4]]);
t3:=Matrix(L,4,4,
[[-zeta4, 0,  0,  0],
 [     0,  1,  0,  0],
 [     0,  0,  0,  1],
 [     0,  0,  1,  0]]);
W2:=MatrixGroup<4,L | [t1,t2,t3]>;

u1:=Matrix(K,4,4,
[ [     1,     0,     0,     0 ],
  [     0,     1,     0,     0 ],
  [     0,     0, zeta3,     0 ],
  [     0,     0,     0,     1 ]]);
u2:=Matrix(K,4,4,
[ [(zeta3+2)/3, (zeta3-1)/3, (zeta3-1)/3, 0 ],
  [(zeta3-1)/3, (zeta3+2)/3, (zeta3-1)/3, 0 ],
  [(zeta3-1)/3, (zeta3-1)/3, (zeta3+2)/3, 0 ],
  [          0,           0,           0, 1 ]]);
u3:=Matrix(K,4,4,
[ [     1,     0,     0,     0 ],
  [     0, zeta3,     0,     0 ],
  [     0,     0,     1,     0 ],
  [     0,     0,     0,     1 ]]);
u4:=Matrix(K,4,4,
[ [(zeta3+2)/3,(1-zeta3)/3, 0,(1-zeta3)/3 ],
  [(1-zeta3)/3,(zeta3+2)/3, 0,(zeta3-1)/3 ],
  [          0,          0, 1,          0 ],
  [(1-zeta3)/3,(zeta3-1)/3, 0,(zeta3+2)/3 ]]);
W3:=MatrixGroup<4,K | [u1,u2,u3,u4]>;

// ProjectiveOrbit computes orbit of 
// points in projective space

ProjectiveOrbit:=function(grp,pt) 
    local i,j,res,v,w,V,PROJ,grpmod,zgrp;
  zgrp:=Centre(grp);
  zgr:=[w : w in zgrp | IsScalar(w)];
  zgrp:=sub<grp | zgrp>;
  grpmod:=Transversal(grp,zgrp); 
  V:=VectorSpace(grp);
  PROJ:=AmbientSpace(Scheme(pt));
  v:=V ! Coordinates(pt);
  res:=[PROJ ! Coordinates(V,v*Transpose(w)) : w in grp];
  return [i : i in Set(res)];
end function;
\end{verbatim}
\end{quotation}

\end{document}